\documentclass[final]{siamart220329}
\usepackage{graphicx}
\usepackage{amsfonts,amsmath,color}
\usepackage{listings}
\usepackage{url,hyperref}
\newcommand{\R}{\mathbb{R}}
\newcommand{\rev}[1]{{\color{black}#1}}

\usepackage[normalem]{ulem}
\newcommand{\bb}[1]{\textcolor{blue}{#1}}

\newcommand{\ignore}[1]{}
\DeclareMathOperator{\rank}{rank}
\DeclareMathOperator*{\minimize}{minimize}

\usepackage{algorithm, algorithmic}
\usepackage{upquote} 

\definecolor{codegreen}{rgb}{0,0.6,0}
\definecolor{codepurple}{rgb}{0.58,0,0.82}
\definecolor{codered}{rgb}{0.7,0,0}

\lstdefinestyle{mystyle}{ 
	commentstyle=\color{codegreen},
	numberstyle=\tiny,
	stringstyle=\color{codepurple},
	basicstyle=\ttfamily\footnotesize,
	breakatwhitespace=false,         
	breaklines=true,                 
	captionpos=b,                   
	keepspaces=true,
	keywordstyle=\color{blue},                                
	showspaces=false,                
	showstringspaces=false,
	showtabs=false,                  
	tabsize=2,
	numbers=left,
	stepnumber=1,
	emph=[1]{feval},emphstyle=[1]\color{black}
}
\title{A Sherman--Morrison--Woodbury approach to solving least squares problems with low-rank updates
}
\author{Stefan G\"uttel\thanks{Department of Mathematics, The University of Manchester, M13 9PL Manchester, United Kingdom, \email{stefan.guettel@manchester.ac.uk},  \email{marcus.webb@manchester.ac.uk}, \email{alban.bloorriley@manchester.ac.uk}. SG is supported by Royal Society Industry Fellowship IF/R1/231032.}\and Yuji Nakatsukasa\thanks{Mathematical Institute, University of Oxford, OX2 6GG Oxford, United Kingdom, \email{nakatsukasa@maths.ox.ac.uk}. Supported by EPSRC grants EP/Y010086/1 and EP/Y030990/1.} \and Marcus~Webb\footnotemark[1]
\and Alban~Bloor~Riley\footnotemark[1] 
}
\date{June 2024}

\begin{document}

\maketitle

\begin{abstract}
We present a simple formula to update the pseudoinverse of a full-rank rectangular matrix that undergoes a low-rank modification, and demonstrate its utility for solving least squares problems. The resulting algorithm can be dramatically faster than solving the modified least squares problem from scratch, just like the speedup enabled by Sherman--Morrison--Woodbury for solving linear systems with low-rank modifications. 
\end{abstract}



\section{Introduction}

The Sherman--Morrison--Woodbury formula (or simply the Woodbury formula)
\[
    (A+UV^T)^{-1}= A^{-1} - A^{-1}U(I + V^T\!A^{-1}U)^{-1}V^T\! A^{-1},
\]
\rev{discovered in the 1950s~\cite{sherman1950adjustment,woodbury1950inverting},} has become a fundamental tool in numerical computation and can be found in many  \rev{popular textbooks; see, e.g., \cite[\S2.1.4]{golub2013matrix}  or \cite[Thm.~18.2.8]{harville1998matrix} or \cite[Thm.~2.3.10]{meyer2023matrix}.} It allows us to efficiently update the inverse of a matrix $A$ when it undergoes a low-rank modification~$A+UV^T$.
Here $A\in\mathbb{R}^{n\times n}$ and $U,V\in\mathbb{R}^{n\times r}$ where usually $r\ll n$. 
The Woodbury formula is particularly useful when $A$ is easy to invert or solve linear systems with\rev{, and often applied for solving linear systems of the form $(A+UV^T)x=b$}.     
It has been used  in numerous applications in scientific computing, including 
quasi-Newton methods~\cite[Ch.~6]{wright1999numerical}, Kalman filtering~\cite{humpherys2012fresh}, and Gaussian processes~\cite{rasmussen2006gaussian}. 
See also  \cite{hager1989updating} for more applications. 

This work was motivated by research into deflation techniques for finding multiple local minima of a nonlinear least squares problem. The proposed deflation operation induces a rank-one \rev{update} to the associated linear least squares problem at each iteration of a Gauss--Newton algorithm~\cite{bloorriley2024deflation}. The Woodbury formula does not apply to the Moore--Penrose pseudoinverse, so cannot be directly applied to these \rev{updated} least squares problems. 

In the 1970s, Meyer~\cite{meyer1973generalized} developed an extension of the Sherman--Morrison formula to the Moore--Penrose pseudoinverse of rectangular matrices when the update is rank one. \ignore{\rev{Wedin~\cite{wedin1973perturbation} 
developed a perturbation theory for pseudoinverses, but the perturbations therein are non-structured (not of low rank) and so the results do not give rise to efficient ways to update solutions of least squares problems.} }
Generalized inverses of square (but possibly singular) matrices updated with blocks of vectors have been discussed in~\cite{henderson1981deriving}.

Despite the prevalence of least squares problems in data science and scientific computing, and despite the existence of Meyer's formula and its variants, to our knowledge pseudoinverse update formulas have not been applied practically to solve least squares problems. In this paper\rev{, using only linear algebra tools accessible to undergraduate students,} we show that such an extension is readily possible, and propose an efficient algorithm WoodburyLS for solving least squares problems wherein  the matrix undergoes a low-rank update. 
The algorithm essentially requires $2r$ solutions of overdetermined least squares problems of the form $\min_x\|Ax-b\|_2$, and 
$2r$ solutions of underdetermined problems of the form $\minimize \|x\|_2$ subject to $A^T x=b$, and in both cases we can reuse a factorization of~$A$ that may have been pre-computed. Alternatively, WoodburyLS can be executed using $r$ solutions of linear systems with the matrix $A^TA$. 
In a typical situation where the QR factorization of $A$ is avilable, the arithmetic cost is reduced by a factor of $O(n/r)$ over the classical solution based on the QR factorization of $A+UV^T$. 
We also 
present a simpler version of Meyer's formula that extends to updates of rank higher than one. We illustrate the performance of the new formula in a simple numerical test, when applied to solving a least squares problem with a low-rank update. 

\ignore{
The Sherman-Morrison formula tells us how the inverse of a matrix gets updated when the matrix inverse is known for example aurora update about diagonal matrix can be inverted in order and squared operations rather than order and cubed. undergoing a low rank update it has been used in various applications including fashion processes optimization other ubiquitous applications. 

Meyer derived the formula, but did not show that it can be used for solving least squares problems efficiently, given the ability to solve problems (least squares and underdetermined problem with $A^T$) with respect to $A$. For example, this is the case when the QR factorization of $A$ has been computed. 
We show that the least squares generalization can be done by solving three \bb{(YN: revise this after converging on algorithm. ) } least squares problems and one  linear system with respect to the transpose of the matrix, where  the minimum-norm solution for an underdetermined linear system is computed. 

In addition, we believe our formula is simpler to understand and implement than Meyer's, who discuss six different cases depending on the properties of the matrix and its relation to the vectors defining the rank-1 update. We derive a simpler formula by focusing on the practically important case where $A$ is tall and full rank. 
}

\section{An update formula for the pseudoinverse}

Given a real\footnote{For simplicity we assume $A,U,V$ are real. For complex matrices, the formulas are valid after replacing the transpose $(\cdot)^T$ with the Hermitian transpose $(\cdot)^*$. 
} matrix $A\in \R^{m\times n}$ 
with $m\geq n$, $\rank(A)=n$, and a vector $b\in\R^m$, we consider the linear least squares problem: find a vector $x\in\R^n$ such that
\[
\| b - Ax \|_2^2 \to \min_x.
\]
It is well known that the 
solution to this problem is given in terms of the Moore--Penrose pseudoinverse by $x = A^\dagger b$~\cite[\S5.5.2]{golub2013matrix}. The standard algorithm is to perform the thin QR factorization $A=QR$ and compute $x=R^{-1}Q^Tb$. 

Let us now consider a low-rank modification of $A$, namely $A + U V^T$ with $U\in \R^{m\times r}$ and $V\in\R^{n\times r}$. We now would like to solve the modified least squares problem, 
\[
\| b - (A + U V^T) \widehat x \|_2^2 \to \min_{\widehat x}.
\]
The vector $b$ could also be different. To address this problem we first present the following theorem, which can be seen as a natural generalization of the Woodbury formula from updates of the matrix inverse to the pseudoinverse.

\begin{theorem}\label{thm:BGNW}
Let $A\in\mathbb{R}^{m\times n}$ and  $U\in\mathbb{R}^{m\times r}, V\in\mathbb{R}^{n\times r}$ with $m\geq n \geq r$. \rev{Assume 
that $\rank(A)=\rank(A+UV^T) = n$.}  
Then
\begin{equation}\label{mainthmmatrix}
(A+UV^T)^\dagger 
= A^\dagger-MA^\dagger
+(I-M)(A^TA)^{-1} VU^T, 
\end{equation}
where
\begin{align*}
    M &= (A^TA)^{-1}  X(I+ Y^T(A^TA)^{-1}X)^{-1} Y^T, \\
    X &= [V, A^TU], \quad\ \ Y = [(A+UV^T)^TU, V].
\end{align*}
\end{theorem}

\smallskip

\ignore{
{\color{red} SG: Just noting that the matrix $I+Y^T (A^T A)^{-1} X$ is invertible iff
$A^T A + X Y^T$ is, the latter of which may be slightly easier to read. Also, $I+Y^T (A^T A)^{-1} X$ 
is the negative Schur complement of the $(2,2)$ block in 
\[
M = \begin{bmatrix}
-I_{2r} & Y^T \\
X &  A^TA
\end{bmatrix}.
\]
Perhaps assuming that $M$ is invertible may lead to a more easy to read condition for our theorem? Can we factor $M$ nicely?
If $\rank(A)=n$, both diagonal blocks in $M$ are invertible. 
Therefore we could use eq.~2 on \url{https://en.wikipedia.org/wiki/Invertible_matrix#Blockwise_inversion} to factor
\[
\begin{bmatrix}
I & -Y^T (A^TA)^{-1}\\
X & I 
\end{bmatrix} M = 
\begin{bmatrix}
I+Y^T (A^T A)^{-1} X & O\\
O & A^T A + X Y^T 
\end{bmatrix}.
\]
}
}

\begin{proof}
Write $\hat A=A+UV^T$. 
As $\hat{A}$ has full rank by assumption, we have 
\begin{align*}
\hat A^\dagger 
&= (\hat A^T\hat A)^{-1} \hat A^T.
\end{align*}
Expanding,
\begin{align*}
\hat A^T\hat A  &= 
A^TA + VU^TA+A^TUV^T
+VU^TUV^T. 
\end{align*}
Now, writing 
$VU^TA+A^TUV^T
+VU^TUV^T = XY^T$ 
where (for example) $X = [V,A^TU]$  
and $Y^T=
\begin{bmatrix}
U^TA+U^TUV^T  \\
V^T
\end{bmatrix}
$, we can apply the Woodbury formula to obtain
\begin{align*}
(\hat A^T\hat A)^{-1}  &= 
(A^TA)^{-1} 
- (A^TA)^{-1}X(I+ Y^T(A^TA)^{-1}X)^{-1}
Y^T(A^TA)^{-1}  \\
&= 
(I-M)(A^TA)^{-1} 
,
\end{align*}
where $M= (A^TA)^{-1}X(I+Y^T(A^TA)^{-1}X)^{-1} Y^T$. 

\rev{Let us convince ourselves that $I + Y^T (A^T A)^{-1} X$ is indeed nonsingular. By definition, $\hat A^T \hat A = A^T A + XY^T$ and so $\rank(\hat A^T \hat A) = \rank(A^T A + XY^T)=n$. Using Sylvester's determinant theorem (see, e.g., \cite[\S18.1]{harville1998matrix}), we have $\det(A^T A + XY^T) = \det(A^TA)\det(I + Y^T (A^T A)^{-1} X)\neq 0$ and hence $I + Y^T (A^T A)^{-1} X$ is indeed nonsingular.

Following on from the above expression for $(\hat A^T\hat A)^{-1}$, we have}
\begin{align*}
\hat A^\dagger
&= (I-M) (A^TA)^{-1} (A+UV^T)^T.
\end{align*}
Finally, we rewrite this in terms of $A^\dagger=(A^TA)^{-1} A^T $:
\begin{align*}
\hat A^\dagger
&= (I-M) (A^TA)^{-1} (A+UV^T)^T
 \\
&= (I-M) (A^TA)^{-1} A^T + (I-M)(A^TA)^{-1} V U^T \\
&= (I-M) A^\dagger + (I-M)(A^TA)^{-1} V U^T.
\end{align*}
This is the update formula for the pseudoinverse given in equation \eqref{mainthmmatrix}.
\end{proof}

Note that~\eqref{mainthmmatrix} is generically a rank-$2r$ update of $A^\dagger$. \rev{The range of the update is contained within the range of $(A^TA)^{-1}X = [(A^TA)^{-1}V, (A^T A)^{-1}A^TU]$, so is of dimension at most $2r$.} By contrast, when $A$ is square the standard Woodbury formula shows that the update to $A^{-1}$ is rank $r$. This is a genuine  difference to keep in mind, but it does not stop us from designing an efficient algorithm for 
solving least squares problems with respect to $\hat A=A+UV^T$ that has a lower computational cost. 

As discussed in the introduction, the topic of low-rank updates of the pseudoinverse has been studied in the literature~\cite{henderson1981deriving,meyer1973generalized}. However, the expressions there appear less suitable for designing an algorithm for least squares problems. 
Note that updating the pseudoinverse of a fat matrix where $m<n$ is simply a matter of transposing the equation~\eqref{mainthmmatrix}. However, cases in which \rev{$A$ or $A+UV^T$ is} rank-deficient are a non-trivial extension we will not discuss in this paper. 
 
\section{Solving updated least squares problems}\label{sec:solvingLSproblems}

To apply Theorem~\ref{thm:BGNW} for the purpose of solving an updated least squares problem $\hat A^\dagger b = (A+UV^T)^\dagger b$, 
we directly use the formula 
\begin{equation}\label{jio}
\hat A^\dagger b = (I - M) (A^\dagger b - (A^T A)^{-1} V U^T b).
\end{equation} 
First note that $M$ involves $(I + Y^T(A^TA)^{-1}X)^{-1}$, which is merely the inverse of a $2r \times 2r$ matrix.

Note also that $M$ involves $(A^TA)^{-1} X = A^\dagger (A^T)^\dagger X$, which requires $2r$ solves with respect to $A^TA$ (one for each column of $X$). At first glance it appears that one would also need to compute $(A^TA)^{-1}V$ directly, but that is unnecessary because $(A^TA)^{-1}V$ can be obtained from the first $r$ columns of $(A^TA)^{-1}X$.

How to do the computation of $(A^TA)^{-1} X$ depends on the situation. For example, if a QR factorization $A=QR$ is available, we can efficiently compute $(A^TA)^{-1} X = R^{-1}R^{-T}X$ via two triangular solves. Otherwise, for example when a preconditioner for $A$ is available for use in an iterative least squares solver (e.g.,~as in \cite{avron2010blendenpik}), we could perform $(A^TA)^{-1} X = A^\dagger (A^T)^\dagger X$ by $2r$ 
solves with respect to $A^T$ and $2r$ solves with respect to $A$ using, e.g.,~LSQR~\cite{paige1982lsqr} and LSRN~\cite{meng2014lsrn}.

To summarize, to solve the least squares problem for $A+UV^T$, we need the solution $x_0$ to the least squares problem for $A$, the solution to a $2r \times 2r$ linear system, and the solution to $2r$ linear systems with the matrix $A^TA$. We provide pseudocode in Algorithm~\ref{pseudocode}, and MATLAB code in Figures~\ref{matlabcode} and~\ref{fig:code}. 

\subsection{Solving for multiple low-rank updates and right-hand sides}

When the right-hand side $b$ stays the same and $U,V$ are modified multiple times, further efficiency savings can be made by precomputing and storing $x_0$ and the machinery for efficiently solving the required least squares systems such as QR factors. 

When one needs to solve least squares problems with $A+UV^T$ for multiple right-hand sides, say $k$ right-hand sides, the cost becomes $2r$ linear systems with $A^TA$ (same as when $k=1$), and the solution of $k$ least squares problems with~$A$. 

In a typical situation where $A$'s QR factorization is given, this means a least squares problem with $A+UV^T$ can be solved in $O((r+k)mn)$ operations instead of the $O(mn^2 + kmn)$ with a standard QR-based approach. In a typical case where $k=O(1)$, this represents a speedup of $O(n/r)$. The complexity of WoodburyLS can be even lower, for example, when $A$ is sparse and well-conditioned so that $A$-solves and $A^T$-solves can be done in $O(\mbox{nnz}(A))$ operations using an iterative solver.

\begin{algorithm}[htbp]
    \caption{WoodburyLS: Solve $\min_x\|b - (A+UV^T)x\|_2$ 
where $A$, $b$, $U$, and $V$ are as in Theorem \ref{thm:BGNW}. Efficient solvers are required for computing $A^\dagger b$ and $(A^T)^\dagger c$, given $b\in\mathbb{R}^{m}$ and $c\in\mathbb{R}^{n}$, or a routine for computing $(A^TA)^{-1}c$.
}\label{alg:WoodburyLS}
  \label{pseudocode}
  \begin{algorithmic}[1] 
    \STATE Compute $x_0=A^\dagger b$, if not already available. 
    \STATE Set $X=[V, A^TU]$, and $Y=[(A+UV^T)U, V]$. 
    \STATE Compute $Z=A^\dagger ((A^T)^\dagger X)$, or alternatively, $Z=(A^TA)^{-1}X$. 
    \STATE Set $w = x_0+Z_1U^Tb$, where $Z_1$ is the first $r$ columns of $Z$. 
    \STATE Compute $\hat w = Z(I_{2r}+Y^TZ)^{-1}Y^Tw$.
    \STATE $x=w-\hat w$ is the solution. 
\end{algorithmic}
\end{algorithm}

\begin{figure}[!ht]
\begin{lstlisting}[language=matlab, numbers=none]
function [x,AtAsolver] = WoodburyLS(A,b,U,V,x0,AtAsolver)
%WoodburyLS Solves the least squares problem min_x ||b-(A+UV')x||
% where A is an m x n matrix, m >= n, and U and V have r columns. 
% Requires A and A+UV' to be full rank.
%
% First call:
%   [x0,AtAsolver] = WoodburyLS(A,b) returns LS solution x0 = A\b
%   and a function AtAsolver that solves A'A x = b for a given b.
%
% Every follow-up call:
%   x = WoodburyLS(A,b,U,V,x0,AtAsolver) returns the LS solution
%   x = (A+U*V')\b using a given solution x0 = A\b and AtAsolver.

if nargin < 3
    [Q,R] = qr(A,0); x = R \ (Q'*b);
    AtAsolver = @(X) R \ (R' \ X);
    return
end
r = size(U,2);
X  = [V, A'*U];                       % X = [V, A'U]
Yt = [X(:,r+1:2*r)' + (U'*U)*V'; V']; % Y = [(A+UV')'U, V]
Z = AtAsolver(X);                     % Z = (A'A)\X
w = x0 + Z(:,1:r) * (U'*b);           % Z(:,1:r) = (A'A)\V
Mw = Z * ((eye(2*r)+Yt*Z) \ (Yt*w));  % M = (A'A)\X(I+Y'(A'A)\X)\Y'
x = w - Mw;
\end{lstlisting}
\caption{MATLAB function implementing 
WoodburyLS described in 
\rev{Algorithm~\ref{alg:WoodburyLS}}.\label{matlabcode}}
\end{figure}



\section{Numerical experiment}

To demonstrate the efficiency gain from using an update formula instead of solving modified problems $(A + U V^T)^\dagger b$ from scratch, we perform a simple experiment as follows.\footnote{Code available at \url{https://github.com/nla-group/WoodburyLS}. \rev{The timings reported here have been produced on a 2023 Windows laptop with an Intel i7-1255U CPU and 32\,GB RAM.}}
Keeping the number of $A$'s rows, $m=10^5$, fixed, we vary both the number of columns $n=100,200,\ldots,1000$ and the rank~$r=10,20,30$ of the update. All matrices $A$ are generated in MATLAB 
with \texttt{randn(m,n)} and \texttt{b=randn(m,1)}. The MATLAB code we used for the timings is essentially given in Figure~\ref{fig:code}, except that we have run each algorithm ten times and averaged the runtimes. We then plot the quotient of the time to compute \texttt{x1} (solving the updated least squares problem from scratch via a QR decomposition of $\hat A = A + U V^T$) over the time required to compute \texttt{x2} (using \texttt{WoodburyLS}). 

\begin{figure}[!ht]
\centering\includegraphics[scale=0.4]{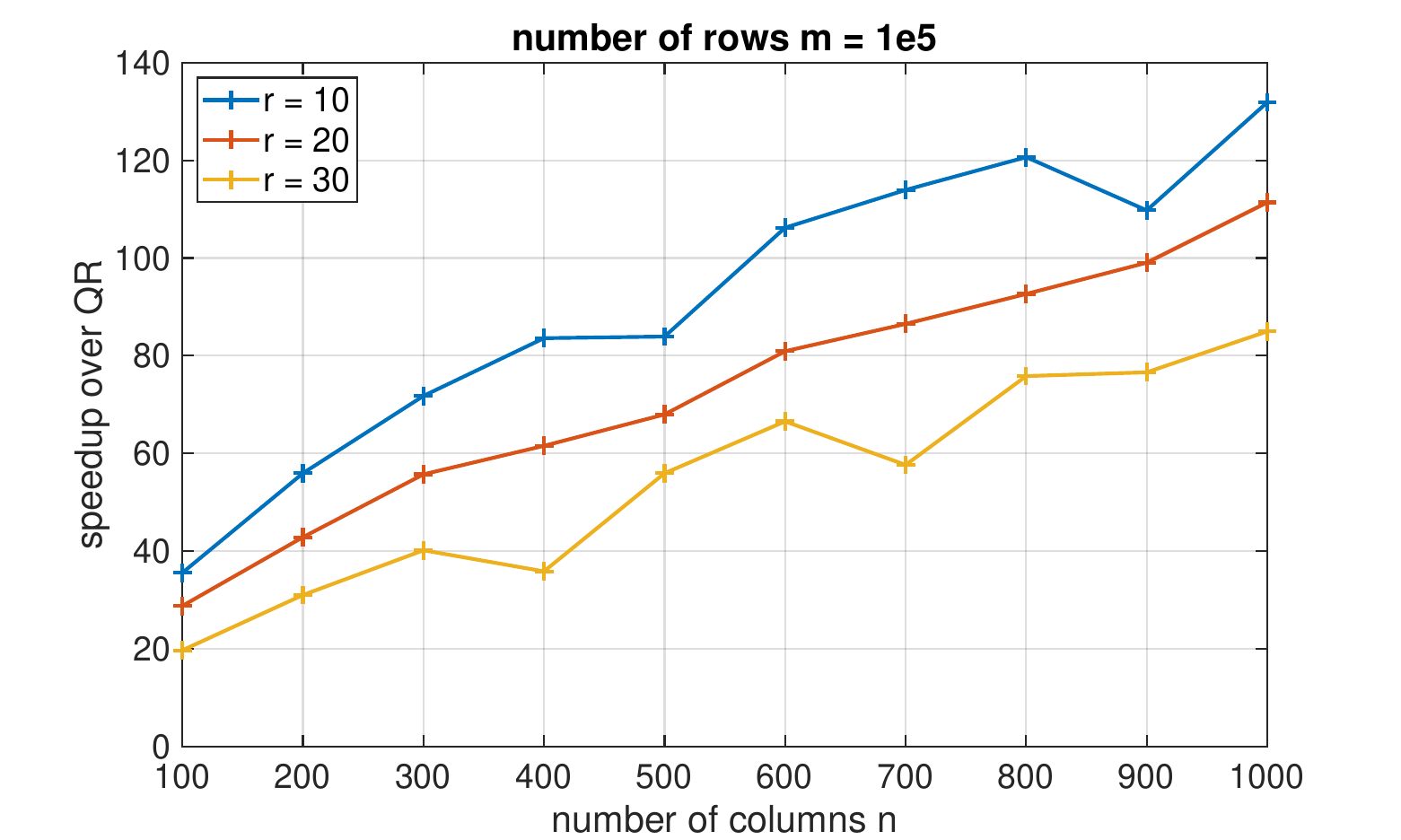}
\caption{Speedup obtained with WoodburyLS compared to solving the modified least squares problem from scratch by computing a QR factorization of $A+U V^T$.}
\label{fig:speedup}
\end{figure}

The results are shown in Figure~\ref{fig:speedup}. In the ranges of parameters tested, we obtain between 20 to over 130-fold speedup. The speedup generally increases with the number of columns~$n$, and it decreases with the rank~$r$. Given that computing \texttt{x1} via QR requires $O(mn^2)$ flops while computing \texttt{x2} via \texttt{WoodburyLS} costs $O(mnr)$, one might expect the speedup in Figure~\ref{fig:speedup} to behave like $O(n/r)$. We find that this is only approximately the case, most likely because flop counts do not directly translate into runtimes due to many other aspects such as memory swaps, blocking and communication costs. \rev{Specifically, the dominant costs are the matrix-matrix product $A^TU$ (costing $\mathcal{O}(nmr)$ operations) and the triangular solves (costing $\mathcal{O}(n^2r)$ operations), which are both highly optimized in MATLAB in the ways just described.} In any case, the practical speedup is significant across the board.
The relative forward error of the computed solution, \texttt{norm(x2-x1)/norm(x1)}, was below $3\times 10^{-14}$ in all cases.

\begin{figure}[!ht]
\begin{lstlisting}[language=matlab, numbers=none]
m = 1e5; n = 500; r = 20;
A = randn(m,n); b = randn(m,1);
U = randn(m,r); V = randn(n,r);

% solve unmodified LS problem via QR
[x0,AtAsolver] = WoodburyLS(A,b);              % 2.340 seconds

% inefficient: min ||b-(A+UV')x|| via QR
Ahat = A + U*V';                               %  
[Qhat,Rhat] = qr(Ahat,0);                      % 2.390 seconds
x1 = Rhat\(Qhat'*b);                           % 

% better solve modified LS problem like this:
x2 = WoodburyLS(A,b,U,V,x0,AtAsolver);         % 0.037 seconds
\end{lstlisting}
\caption{\label{fig:code} An example demonstrating the use of \texttt{WoodburyLS} in MATLAB. Both \texttt{x1} and \texttt{x2} are solutions to the modified least squares problem $\min_x \| b - (A+UV^T)x \|_2$, but \texttt{x2} is computed significantly faster.}
\end{figure}

\section*{Acknowledgements}
We are grateful to the organisers of the British Mathematical Colloquium 2024, held at The University of Manchester in June 2024, where much of the research for this paper took place.
\ignore{\rev{We thank Ilse Ipsen for bringing \cite{wedin1973perturbation} to our attention.}}

\bibliographystyle{siamplain}
\bibliography{main}

\begin{thebibliography}{10}

\bibitem{avron2010blendenpik}
{\sc H.~Avron, P.~Maymounkov, and S.~Toledo}, {\em Blendenpik: Supercharging
  {LAPACK}'s least-squares solver}, SIAM J. Sci. Comp., 32 (2010),
  pp.~1217--1236.

\bibitem{bloorriley2024deflation}
{\sc A.~Bloor~Riley, M.~Webb, and M.~Baker}, {\em Deflation techniques for
  finding multiple local minima of a nonlinear least squares problem}, in
  prep.,  (2024).

\bibitem{golub2013matrix}
{\sc G.~H. Golub and C.~F. Van~Loan}, {\em Matrix Computations, 4th Edition},
  John Hopkins University Press, 2013.

\bibitem{hager1989updating}
{\sc W.~W. Hager}, {\em Updating the inverse of a matrix}, SIAM Rev., 31
  (1989), pp.~221--239.

\bibitem{harville1998matrix}
{\sc D.~A. Harville}, {\em Matrix Algebra from a Statistician's Perspective},
  Taylor \& Francis, 1998.

\bibitem{henderson1981deriving}
{\sc H.~V. Henderson and S.~R. Searle}, {\em On deriving the inverse of a sum
  of matrices}, SIAM Rev., 23 (1981), pp.~53--60.

\bibitem{humpherys2012fresh}
{\sc J.~Humpherys, P.~Redd, and J.~West}, {\em A fresh look at the kalman
  filter}, SIAM Rev., 54 (2012), pp.~801--823.

\bibitem{meng2014lsrn}
{\sc X.~Meng, M.~A. Saunders, and M.~W. Mahoney}, {\em {LSRN}: A parallel
  iterative solver for strongly over-or underdetermined systems}, SIAM J. Sci.
  Comp., 36 (2014), pp.~C95--C118.

\bibitem{meyer1973generalized}
{\sc C.~D. Meyer}, {\em Generalized inversion of modified matrices}, SIAM J.
  Appl. Math., 24 (1973), pp.~315--323.

\bibitem{meyer2023matrix}
{\sc C.~D. Meyer}, {\em Matrix Analysis and Applied Linear Algebra, 2nd
  Edition}, SIAM, 2023.

\bibitem{wright1999numerical}
{\sc J.~Nocedal and S.~J. Wright}, {\em Numerical Optimization}, Springer New
  York, second~ed., 1999.

\bibitem{paige1982lsqr}
{\sc C.~C. Paige and M.~A. Saunders}, {\em {LSQR}: An algorithm for sparse
  linear equations and sparse least squares}, ACM Trans. Math. Soft., 8 (1982),
  pp.~43--71.

\bibitem{rasmussen2006gaussian}
{\sc C.~E. Rasmussen and C.~K.~I. Williams}, {\em Gaussian Processes for
  Machine Learning}, MIT Press Cambridge, 2006.

\bibitem{sherman1950adjustment}
{\sc J.~Sherman and W.~J. Morrison}, {\em Adjustment of an inverse matrix
  corresponding to a change in one element of a given matrix}, Ann. Math.
  Stat., 21 (1950), pp.~124--127.

\bibitem{woodbury1950inverting}
{\sc M.~A. Woodbury}, {\em Inverting modified matrices}, Stat. Res. Group Memo
  Repos., Princeton University, 42 (1950), pp.~1--4.

\end{thebibliography}

\end{document}